\documentclass{amsart}

\issueinfo{00}% volume number
  {}%         % issue number
  {}%         % month
  {1995}%     % year

\theoremstyle{definition}

\theoremstyle{remark}

\numberwithin{equation}{section}

\begin{document}

\title{White Hole, Black Whole, and The Book}

% author one information
\author{K. K. Nambiar}
\address{Formerly, Jawaharlal Nehru University, New Delhi, 110067, India}
\curraddr{1812 Rockybranch Pass, Marietta, GA, 30066-8015, U.S.A.}
\email{nambiar@mediaone.net}
%\thanks{}

% author two information
%\author{}
%\address{}
%\email{}
%\thanks{}

\subjclass{Primary 03A05, 03E30; Secondary 03E17, 03E50}
%\subjclass{Primary 68P10, 68R10; Secondary 68P30, 15A48}
\date{February 18, 2001}

\dedicatory{Dedicated to the memory of Professor 
Paul Erd\"os, the originator of The Book.}

% presently the "communicated by" line appears only in PROC
%\commby{}

\begin{abstract}
Physical and intellectual spaces are visualized 
making use of concepts from intuitive set theory. Intellectual 
space is defined as the set of all proofs of mathematical 
logic, contained in The Book conceived by Erd\"os.
\vskip 5pt \noindent
{\it Keywords\/}---Physical space, Intellectual space, Visualization.
\end{abstract}

\maketitle

\specialsection{INTRODUCTION}

In an earlier paper \cite {Nam:IST}, it was shown that 
Zermelo-Fraenkel set theory 
gets considerably simplified, if we add two axioms, 
\emph {Monotonicity} and \emph {Fusion}, to it. In the resulting 
intuitive set theory (IST), the continuum hypothesis is a theorem, 
axiom of choice is a theorem, Skolem Paradox does not crop up,
non-Lebesgue measurable sets are not possible, and the unit interval 
splits into a set of infinitesimals with cardinality 
\(\aleph_0\) \cite {Nam:IST,Nam:VIST}. 
This paper shows that IST can be used to 
visualize the infinite physical space around us as a set. Further, 
if we consider all the proofs of mathematics as our intellectual 
space, then IST provides a way to consider that also as a set.

\specialsection{WHITE HOLE}

The axiom of fusion allows us to imagine a unit interval as a set 
of infinitesimals, with each infinitesimal 
containing \(\aleph_{1}\) 
\emph {figments} (elements which cannot be accessed by the axiom of 
choice) in it. We consider these infinitesimals as integral units 
which cannot be broken up any further. To facilitate the discussion,
in addition to Dedekind 
cuts, we will use also the concept of a \emph {Dedekind knife}, and 
assume that the knife can cut any interval given to it, 
\emph {exactly} in the middle. From this, it follows that every infinite 
\emph {recursive}
subset of positive integers, or equivalently, a binary number in 
the unit interval, represents the use of Dedekind knife an infinite 
number of times. The result we get when we use the knife $\aleph_0$ 
times, according to an infinite binary sequence, is what 
we call an infinitesimal and the location of that infinitesimal is 
what we call a number. 

What if, the operation of the knife is 
continued further an infinite number of times 
according to a new arbitrary infinite binary sequence. We can see the 
intuitionists protesting at this stage, that you cannot start another 
infinite sequence before, you have \emph {completed} the previous infinite 
sequence. For this, the formalist answer is that, in mathematics, 
there is no harm in imagining things which cannot be accomplished 
physically. We can see here, the source of the oxymoron \emph 
{completed infinity}, and the motivation for our definition of a 
\emph {bonded set}. Bonded set is a set, from which axiom of choice cannot 
pick an element and separate it. These special elements, we call 
figments.

Having stated this, we continue our second infinite cutting, this 
time, without bothering to restrict ourselves to recursive subsets 
of positive integers as in the original case.
The justification for this is that our operation is in the realm of the 
imaginable and not physical. The result of the cutting is $\aleph_1$ 
figments and they are to be considered 
only as a figment of our imagination. These arguments, 
of course, do not prove the consistency of the axiom of fusion, but 
hopefully makes it plausible.

The discussion above allows us to define a \emph {white hole} 
(WhiteHole, whitehole)
as the infinitesimal (bonded set) corresponding to an infinite recursive
subset of positive integers (a binary number in the unit interval).
It represents an indefinitely small void, which cannot be broken up 
any further. However, it does contain $\aleph_1$ figments which cannot 
be isolated.

\specialsection{BLACK WHOLE}

A binary number is usually defined as a two way infinite binary sequence 
around the binary point, 
\[
\ldots 000xx\ldots xxx.xxxxx\ldots
\]
in which the $x$s represent either a $0$ or $1$, and the infinite 
sequence on the left eventually ends up in $0$s. The two's 
complement number system represents a
negative number by a two way infinite sequence,
\[
\ldots 111xx\ldots xxx.xxxxx\ldots
\]
in which the infinite sequence on the left eventually ends up in $1$s.

We define the \emph {Universal Number System} as the number system 
in which there are no restrictions on the infinite 
sequences on both sides.

It is easy to recognize that, the sequence 
\[
\ldots 00000.xxxxx\ldots
\]
with a nonterminating binary sequence on the right side
represents a number in a unit interval and also a whitehole. The concept 
of a \emph {black whole} (BlackWhole, blackwhole) is now easy to define.
The two way infinite sequence we get when we flip the whitehole around 
the binary point,
\[
\ldots xxxxx.00000\ldots 
\]
represents a \emph {supernatural number} and also a \emph {black stretch}.
The infinite set of black stretches, we define as the blackwhole. Thus, 
blackwhole can be considered as a dual of the unit interval. The name 
black stretch is supposed to suggest that it can be visualized as a set
of points distributed over an infinite line, but it should be recognized
that it is a bonded set, which the axiom of choice cannot access. Our 
description of the black whole clearly indicates that it can be used 
to visualize what is beyond the finite physical space around us.

\specialsection{THE BOOK}

\emph {The Book} (TheBook, the book) as originally conceived by Erd\"os is 
a book that contains all 
the smallest proofs of mathematics arranged in the lexical order. Since 
the alphabet of any axiomatic theory is finite, and every proof is a 
well-defined formula, it follows that a computer can be set up
to start writing this book. We cannot expect the computer to stop, since 
there are an infinite number of proofs in mathematics. Thus, a computer 
generated book will always have to be unfinished, the big difference 
between a computer generated book and The Book is that it is a 
\emph {finished} book. 

The physical appearance of the book can be visualized as below. 

\begin{itemize} 
\item The front cover and the back cover are each one millimeter 
thick, and the entire book, including the covers, is three millimeters 
thick.
\item The first sheet of paper is half-millimeter thick, the 
second sheet is half thick as the first, the third sheet is 
half thick as the second, and so on.
\item On every odd page is written a full proof, and in the next
even page is written the corresponding theorem. 
\item The last sheet is stuck with the cover with the result that 
the Last Theorem is not visible.
\end{itemize}

From the description of the book, we can infer that any formula 
which is a theorem can be found in the book, by sequentially
going through the pages of the book. The only difficulty 
is that, if a formula 
is not a theorem, we will be eternally searching for it. 

\specialsection{CONCLUSION}

The upshot of all our discussion, reduces to the following. We live
in a space where the visible finite part is filled up with white 
holes. When we use Dedekind knife 
an infinite number of times 
to cut a line segment 
we get a real number and a white hole. 
The invisible part of the physical 
space is an unimaginably large, unreachable black whole, comprising
of transfinite black stretches. The Book allows us to read through 
proofs and collect as many theorems as we want, but it does not help 
us very much in deciding whether a given formula is a theorem. The
main problem of mathematics is to write a New Book with the \emph {theorems} 
listed in lexical order. Hilbert once had hopes of setting up a computer 
to start writing this book, but a great achievement of the twentieth 
century is the dashing of that hope.

%\bibliography{references}
\bibliographystyle{amsplain}

\end{document}